\documentclass[11pt]{amsart}
\usepackage{times}

\newtheorem*{Thm}{Theorem}
\newtheorem*{cor}{Corollary}

\theoremstyle{definition} 
\newtheorem{Def}{Definition}
\newtheorem{lem}[Def]{Lemma}
\newtheorem{prop}[Def]{Proposition}
\newtheorem{rem}[Def]{Remark}

\theoremstyle{remark} 

\newtheorem*{ack}{Acknowledgment}

\def \C {\mathbb{C}}
\def \R {\mathbb{R}}
\def \Q {\mathbb{Q}}
\def \Z {\mathbb{Z}}

\def \F {\mathbb{F}}
\def \a {\alpha}
\def \b {\beta}

\def \del {\partial}

\newcommand{\inner}[2]{\langle {#1},{#2} \rangle}
\newcommand{\exact}[1]{\stackrel{#1}{\longrightarrow}}
\newcommand{\dlim}{\displaystyle\lim_{\longrightarrow}}

\DeclareMathOperator{\tr}{Tr}

\DeclareMathOperator{\Int}{Int}
\DeclareMathOperator{\id}{Id}

\begin{document}

\title[A functional equation for the Lefschetz zeta functions]{A functional equation for the Lefschetz zeta functions of infinite cyclic coverings with an application to knot theory}

\author{Akio Noguchi}
\address{
Department of Mathematics,
Tokyo Institute of Technology,
Oh-okayama, Meguro, Tokyo 152-8551, Japan
}
\email{akio@math.titech.ac.jp }
\thanks{The author was supported in part by JSPS fellowship for young scientists.}


\subjclass[2000]{Primary 57M27; Secondary 37C30}

\begin{abstract}
The Weil conjecture is a delightful theorem for algebraic varieties on finite fields and an important model for dynamical zeta functions. In this paper, we prove a functional equation of Lefschetz zeta functions for infinite cyclic coverings which is analogous to the Weil conjecture. Applying this functional equation to knot theory, we obtain a new view point on the reciprocity of the Alexander polynomial of a knot.
\end{abstract}

\maketitle

\section*{Introduction}\label{sec:intro}
The Lefschetz zeta function is one of the dynamical zeta functions. Dynamical zeta functions are developed in order to study the number of fixed points or periodic points. These zeta functions are motivated by the Weil conjecture. Therefore, studying dynamical zeta functions is usually modeled after Weil conjecture.

On the other hand, dynamical zeta functions are useful in geometric topology. For example, if you fix a map to a geometrical one, it can be related to topological invariants (e.g. \cite{MR2001a:37031,MR2000h:57063,MR2000i:57055,MR2003c:57032}). That is, a property of a topological invariant can be regarded as a property of the dynamical zeta function.

Grothendieck (after the works of Weil and Serre) observed that the Weil conjecture should be the consequence of a certain good cohomology theory, which is called the Weil cohomology, and in particular the functional equation should be derived from the Poincar\'e duality of the cohomology (see section 1.2). In this paper, we prove an analogous functional equation of Lefschetz zeta functions by following his idea, and apply it to knot theory. To do that, we have to resolve the following two problems: 
\begin{enumerate}
\item the infinite cyclic covering of a knot compliment is a non-compact odd-dimensional manifold but the Weil cohomology requires the even-dimensional Poincar\'e duality, and
\item this manifold has the boundary, so we need to deal with the relative version of Poincar\'e duality (Lefschetz-Poincar\'e duality).
\end{enumerate}

The answer of the first is Milnor's duality theorem for infinite cyclic coverings \cite{MR39:3497}. The device for the second is the Lefschetz zeta function for the boundary, which is defined in Definition~\ref{rel-res}. With those tools, the following theorem is proved in Section~\ref{sec:equ}.

\begin{Thm}
Let $M$ be a compact connected manifold of dimension $n$ which may have a boundary $\partial{M}$,
and $\tilde{M}$ an orientable infinite cyclic covering of $M$ with $\dim H_* (\tilde{M};\Q) < \infty$. 
Let $f:(\tilde{M},\partial\tilde{M}) \to (\tilde{M},\partial\tilde{M})$
be a proper continuous map with degree $\lambda \ne 0$ with respect to the compact support cohomology $H^n_{\rm{cpt}} (\tilde{M},\partial \tilde{M};\Q)$. If $n$ is odd, then the Lefschetz zeta function $\zeta_f$ and $\zeta_{f|\partial\tilde{M}}$ satisfies the following functional equation:
$$
\frac{{\zeta_f({1}/{\lambda z})}^2 }{ \zeta_{f|\partial \tilde{M}}({1}/{\lambda z})}
	= {\lambda}^{\chi} z^{2\chi} \frac{{\zeta_f(z)}^2 }{ \zeta_{f|\partial \tilde{M}}(z)},
$$
where $\chi$ is the Euler characteristic of $\tilde{M}$.
\end{Thm}

In the Main theorem, the Lefschetz zeta function for the boundary is essential. However, it is usually easy to compute it. In section 3, we apply the functional equation to knot theory. Computing directly the Lefschetz zeta function for the boundary, we obtain the following well-known theorem~\cite{seifert}.

\begin{cor}
The Alexander polynomial $\Delta_K(z)$ of a knot $K$ is reciprocal.
\end{cor}

This corollary implies that the reciprocity of the Alexander polynomial of a knot is a special case of the functional equation of Lefschetz zeta functions for infinite cyclic coverings. Formally, the reciprocity of the Alexander polynomial was interpreted as symmetry of the covering transformations. However, not only covering transformations but also more general maps have symmetry, and those symmetries are realized in the functional equation. That is, this functional equation would be useful in studying the symmetry of a knot.

The covering transformation of the infinite cyclic coverings plays an analogous role to the Frobenius automorphism in this paper. The similarity between the Lefschetz zeta function and the Reidemeister torsion (also the Alexander polynomial) was first pointed out by Milnor~\cite{MR39:3497} and studied by some other investigators from many viewpoint. For example, Franks~\cite{MR83h:58074} studied the Lefschetz zeta function of a flow on $S^3$ which is associated with a knot, and as a result he reproved the reciprocity of the Alexander polynomial. We also refer to Morishita's article \cite{preMori} for some other analogies with Iwasawa theory.

Grothendieck also studied $p$-adic cohomology theory to approach the Weil conjecture. That immediately asks whether $p$-adic (co)homology theory is also helpful for knot theory. The answer is given in \cite{math.GT/0503252}. The $p$-adic coefficient Alexander module is helpful to investigate what the zeros of $\Delta_K(z)$ mean.

\bigskip

The first section is a preliminary section. In section~\ref{sec:dynam}, we recall the Lefschetz zeta function. In section~\ref{sec:weil}, we review the Weil conjecture. In section~\ref{sec:duality}, we recall a duality theorem for infinite cyclic coverings of manifolds. This theorem was obtained by Milnor~\cite{MR39:3497}. The proof of Main Theorem is in section~\ref{sec:equ}. We prove this theorem, using the duality of infinite cyclic coverings in place of the Weil cohomology's duality. In section~\ref{sec:knot}, we apply the functional equation of the main theorem to knot theory.

\section{preliminaries}
\subsection{Dynamical zeta functions}\label{sec:dynam}
Dynamical zeta functions are types of zeta function which have been developed in order to study the number of fixed points and periodic points of a map. The most fundamental dynamical zeta function is defined by Artin and Mazur~\cite{MR31:754}. This zeta function is called the Artin-Mazur zeta function, and contains all the information about the numbers of fixed points for iteration.
While the Artin-Mazur zeta function counts the periodic points geometrically, the Lefschetz zeta function counts them homologically.

\begin{Def}\label{Lef}~\nocite{MR37:3598}
Let $M$ be a manifold with $\dim H_* (M,\Q) < \infty$,
and $f:M \to M$ a continuous map on $M$.
The {\it Lefschetz zeta function} is defined by
$$
\zeta_f(z)=\exp\sum_{n=1}^{\infty} \frac{\Lambda(f^n)}{n} z^n ,
$$
where
$$\Lambda(f)=\sum_{i=0}^{\dim M} (-1)^i \tr (f_{*i}:H_i(M;\Q) \to H_i(M;\Q)) $$
is the Lefschetz number of $f$.
\end{Def}

The Lefschetz zeta function has some better properties than the Artin-Mazur zeta function. Not all Artin-Mazur zeta functions have a strictly positive radius of convergence, but Lefschetz zeta functions always have. This means that the Artin-Mazur zeta function is a `formal function', but the Lefschetz zeta function is always a function. Moreover, the Lefschetz zeta function is always a rational function. (Artin-Mazur zeta functions for Axiom A diffeomorphisms are rational functions~\cite{MR44:5982}.) This rationality is a remarkable property so that Lefschetz zeta functions can be related to characteristic polynomials. (The proof is written in~\cite{MR37:3598} for example.)

\begin{prop}\label{rationality}
The Lefschetz zeta function admits a rational expression
$$
\zeta_f(z)=\prod_{i=0}^{\dim X} \det(\id - z f_{*i})^{(-1)^{i+1}},
$$
where $\id$ is defined as the identity map. In particular, it is always a rational function.\end{prop}

\begin{rem}\label{rationality2}
Since we use homology with rational coefficients, homology groups can be replaced by cohomologies.
\end{rem}

\subsection{Weil conjecture}\label{sec:weil}

Our theorem is an analogy to the Weil conjecture.
In this section we recall the Weil conjecture and related topics briefly.

The congruence zeta function arises from the study of the number of solutions of a congruence
$$
f(x_1,\dots,x_n) \equiv 0 \quad \mod p,
$$
where $p$ is a prime and $f(x_1,\dots,x_n)$ is a polynomial of integral coefficients. It is natural to look for solutions not only in the prime field $\F_{p}$ but also in all of its finite extensions $\F_{p^m}$.
Abstracting this idea to algebraic varieties, the congruence zeta function of algebraic variety is defined as follows.

\begin{Def}
Let $\F_{p}$ be a finite field with $q$ elements and $V$ a algebraic variety of dimension $n$ defined over $\F_{p}$. Let $\F_{p^m}$ be the extension field of $\F_{p}$ of degree $m$ and $N_m$ the number of $\F_{p^m}$-rational points of $V$. 
Then the function $Z(u,V)$ of $u$ defined by
$$
Z(u,V)=\exp \sum_{m=1}^\infty \frac{N_m}{m} u^m
$$
is called the {\it congruence zeta function} of the algebraic variety $V$.
\end{Def}

The congruence zeta function $Z(u,V)$ has the following properties. These properties were conjectured by Weil~\cite{MR10:592e} and finally proved by Deligne~\cite{MR49:5013}. But this theorem is still called the Weil conjecture.

\begin{prop}[Weil conjecture]\label{Weil}
If $V$ is a nonsingular projective variety over $k$, then $Z(u,V)$ has the following properties:\\
(1) {\it Rationality}:
$Z(u,V)$ is a rational function in $u$.\\
(2) {\it Functional equation}:
$Z(u,V)$ satisfies the functional equation
$$
Z(\frac{1}{q^n u},V)=\pm q^{n\chi/2} u^\chi Z(u,V),
$$
where the integer $\chi$ is the Euler characteristic of V.\\
(3) {\it Riemann hypothesis}:
$Z(u,V)$ can be factored as
$$
Z(u,V)=\prod^{2n}_{i=0} P_i(u)^{(-1)^{i+1}},
$$
where $P_i(u)=\prod_{j=1}^{B_i} (1-\alpha_j^{(i)} u)$ and $\alpha_j^{(i)}$ is satisfies $|\alpha_j^{(i)}|=q^{i/2}$.
\end{prop}

Now let us review how the Weil conjecture was solved. Weil initially suggested the possibility of usung the Lefschetz fixed point formula to approach the Weil conjecture~\cite{MR0092196}. It means that the congruence zeta functions can be interpreted as the Lefschetz zeta function of the Frobenius action. (In view of this, he should get credit for the definition of the Lefschetz zeta function.)

The next, inspired by works of Serre, Grothendieck formulated the cohomology theory which is required to realize the Weil conjecture, which is called Weil cohomology (see below), and studied the \'etale cohomology toward realizing the Weil cohomology~\cite{SGA}. Actually, his \'etale cohomology proved the functional equation. 

\begin{Def}[Weil cohomology]
Let $K$ be a field of characteristic 0. A contravariant functor $V \to H^*(V)$ is called a {\it Weil cohomology} with coefficients in $K$ if it has the following three properties:

(1) {\it Poincar\'e duality}: If $n=\dim V$, then a orientation isomorphism $H^{2n}(V) \cong  K$ exits and the cup product $H^j(V) \times H^{2n-j}(V) \to H^{2n}(V) \cong K$ induces a non-degenerate pairing.

(2) {\it K\"unneth formula}: For any $V_1$ and $V_2$ the mapping $H^*(V_1) \otimes H^*(V_2) \to H^*(V_1 \times V_2)$ defined by $a \otimes b \mapsto Proj_1^*(a) \cdot Proj_2^*(b)$ is an isomorphism.

(3) {\it Cycle map}: Let $C^j(V)$ be the group of algebraic cycles of codimension $j$ on V. There exists a fundamental class homomorphism $FUND: C^j(V) \to H^{2j}(V)$ for all $j$ which is functorial, compatible with products via the K\"unneth map, and sends a 0-cycle to its degree as an element of $H^{2n}(V)$.
\end{Def}

The Weil cohomology was modeled after the classical cohomology theory of smooth compact varieties over the complex numbers $\C$. Actually, in this case the cohomology theory satisfies the condition of the Weil cohomology.

The expository articles~\cite{MR45:1920,MR51:8113} are helpful to see how the Weil cohomology derives the Weil conjecture. The point is that the even-dimensional Poincar\'e duality derives the functional equation. (In view of this, W\'ojcik's functional equation in~\cite{MR97i:58141} might have been already known in Grothendieck's project. However his computation is helpful in this paper.)

\subsection{Duality theorem for infinite cyclic coverings}\label{sec:duality}
To prove functional equations, we need an even-dimensional Poincar\'e duality. But the infinite cyclic covering of a knot compliment is three dimensional and apparently does not admit the functional equation. However Milnor~\cite{MR39:3497} proved a duality theorem, which look like even-dimensional Poincar\`e duality for odd-dimensional manifolds. In this section, we give a brief review of this duality.

Let $M$ be a manifold. The {\it infinite cyclic covering} of $M$ is defined as covering space $\tilde{X}$ which is determined by some homomorphism of the fundamental group $\pi_1 (X)$ onto an infinite cyclic group. These spaces have two ends.

Let $N_\a$ and $N'_\a$ be neighborhoods of two ends $\epsilon$ and $\epsilon'$, and $\{N_\a \cup N'_\a\}$ be directed. The direct limit of the Mayer-Vietris sequence:
\begin{equation*}
\begin{split}
\cdots \exact{} H^{i-1}(\tilde{M},N_\a \cap N'_\a) \exact{\delta^*} H^i&(\tilde{M},N_\a \cup N'_\a) \exact{} \\& H^i(\tilde{M},N_\a) \oplus H^i(\tilde{M},N'_\a) \exact{} \cdots,
\end{split}
\end{equation*}
prove that the connecting homomorphism
$$
\delta^*(=\dlim \delta^*): H^{i-1}(\tilde{M}) \to H^i_{\rm{cpt}}(\tilde{M})
$$
is an isomorphism. This isomorphism is essential.

From the Poincar\'e duality theorem for an oriented $n$-dimensional manifold, the cup product
$$
\cup: H^i_{\rm{cpt}} (\tilde{M}) \times H^{n-i} (\tilde{M},\partial \tilde{M}) \to H^{n}_{\rm{cpt}} (\tilde{M},\partial \tilde{M}) \cong \Q
$$
provides a non-degenerate pairing. 
By the above argument, we have the following.

\begin{prop}[Duality theorem for infinite cyclic coverings \cite{MR39:3497}]\label{duality}
Let $M$ be a compact connected $n$-dimensional manifold with boundary,
and $\tilde{M}$ an orientable infinite cyclic covering of $M$.
If $H_*(\tilde{M};\Q)$ is finitely generated over $\Q$,
then the cup product:
$$
\cup: H^{i-1} (\tilde{M};\Q) \times H^{n-i}(\tilde{M},\partial \tilde{M};\Q) \to H^{n-1}(\tilde{M},\partial \tilde{M};\Q)\cong \Q
$$
provides a non-degenerate pairing.
\end{prop}

\begin{rem}
In original paper, $M$ is assumed to be triangulated. However, from the study of Kirby and Siebenmann~\cite{MR58:31082}, it is not necessary. It was pointed out in \cite{MR58:7645}\cite{MR1102252}.
\end{rem}

\section{Functional equation for infinite cyclic coverings}\label{sec:equ}
In this section, we prove our main theorem, replacing the the Poincar\`e duality of the Weil cohomology with Milnor's duality theorem (Proposition~\ref{duality}). However, only one problem remains. We need to deal with the relative version of Milnor duality because the infinite cyclic covering has the boundary. For that reason, we introduce other Lefschetz zeta functions in this section (Definition~\ref{rel-res}). After that, we prove the following theorem.

\begin{Thm}
Let $M$ be a compact connected manifold of dimension $n$ which may have a boundary $\partial{M}$, and $\tilde{M}$ an orientable infinite cyclic covering of $M$ with $\dim H_* (\tilde{M};\Q) < \infty$.
Let $f:(\tilde{M},\partial\tilde{M}) \to (\tilde{M},\partial\tilde{M})$
be a proper continuous map with degree $\lambda \ne 0$ with respect to the compact support cohomology $H^n_{\rm{cpt}} (\tilde{M},\partial \tilde{M};\Q)$. If $n$ is odd, then the Lefschetz zeta function $\zeta_f$ and $\zeta_{f|\partial\tilde{M}}$ satisfies the following functional equation:
$$
\frac{{\zeta_f({1}/{\lambda z})}^2 }{ \zeta_{f|\partial \tilde{M}}({1}/{\lambda z})}
	= {\lambda}^{\chi} z^{2\chi} \frac{{\zeta_f(z)}^2 }{ \zeta_{f|\partial \tilde{M}}(z)},
$$
where $\zeta_{f|\partial \tilde{M}}(z)$ is the restricted Lefschetz zeta function to the boundary (cf. Definition~\ref{rel-res}), and $\chi$ is the Euler characteristic of $\tilde{M}$.

In particular, in the case where $\partial \tilde{M} = \emptyset$, we have
$$
{\zeta_f(\frac{1}{\lambda z})}
	= \pm {\lambda}^{\chi/2} z^{\chi} {\zeta_f(z)}
$$(cf.~Proposition~\ref{Weil}).
\end{Thm}

We define two other Lefschetz zeta functions (cf. Definition~\ref{Lef}). These zeta functions let us use the property of the Lefschetz numbers in the proof of the Main Theorem.

\begin{Def}\label{rel-res}
Let $(M,A)$ be a pair of manifolds with $\dim H_* (M,A;\Q) < \infty$ and $\dim H_* (A;\Q) < \infty$.
Let $f:(M,A) \to (M,A)$ be a continuous map.
The {\it relative Lefschetz zeta function} is defined by
$$
\zeta_f^{\rm{rel}}(z) = \exp\sum_{n=1}^{\infty} \frac{\Lambda^{\rm{rel}}(f^n)}{n} z^n ,
$$
where $\Lambda^{\rm{rel}}(f) = \sum_{i=0}^{\dim M} (-1)^i \tr(f_{*i}:H_i(M,A;\Q) \to H_i(M,A;\Q))$.
The {\it restricted Lefschetz zeta function to $A$} is defined by
$$
\zeta_{f|A}(z) = \exp\sum_{n=1}^{\infty} \frac{\Lambda((f|A)^n)}{n} z^n ,
$$
where $\Lambda(f|A) = \sum_{i=0}^{\dim A} (-1)^i \tr((f|A)_{*i}:H_i(A;\Q) \to H_i(A;\Q)).$
\end{Def}

\begin{lem}\label{form:rel}
If two of the three homology groups $H_{*} (M;\Q), H_{*} (A;\Q)$ and $H_{*} (M,A;\Q)$ are finite dimensional, then all three Lefschetz zeta functions are defined and satisfy the following relation:
$$\zeta_f (z) = \zeta^{\rm{rel}}_f (z) \times \zeta_{f|A} (z).$$
\end{lem}

\begin{proof}
We can see that $\Lambda(f)=\Lambda^{\rm{rel}}(f)+\Lambda(f|A)$, which follows from the exact sequence:
$$
\cdots \exact{\del_*}  H_{i} (A) \exact{\a_*}  H_{i} (X) \exact{\b_*}  H_{i} (X,A) \exact{\del_*} \cdots.
$$
From the definitions, we can see that
$\zeta_f (z) = \zeta^{\rm{rel}}_f (z) \times \zeta_{f|A} (z)$.
\end{proof}

We need some formulas to prove the main theorem.

\begin{lem}[\cite{MR97i:58141}]\label{formula}
Let $\inner{~}{~}:V \times V' \to \Q$ be a non-degenerate pairing of the $\Q$-vector spaces with dimension $n$.
Let $\lambda \in \Q \setminus \{0\}$, and $f:V \to V, g:V' \to V'$ be endmorphisms such that
$$
\inner{f(x)}{g(y)} = \lambda \inner{x}{y}
$$
for all $x \in V,y \in V'$. Then
\begin{align*}
&(1)& &(\det f)(\det g)=\lambda^n,\\
&(2)& \det(\id - gt)&\det f =(-1)^n \lambda^n t^n \det(\id -f/\lambda t).
\end{align*}
Here $t$ is an indeterminacy.
\end{lem}

\begin{proof}
Since $\lambda \in \Q \setminus \{0\}$, $f$ and $g$ are isomorphisms, and their inverse maps $f^{-1},g^{-1}$ exist.
The equation
$
\inner{x}{\lambda g^{-1}(y)} = 
\inner{f(x)}{y}
$
means that the dual map $f^* = \lambda g^{-1}$. 

Hence we obtain 
\begin{align*}
\det f = \det f^* = \det \lambda g^{-1} = \lambda^n / \det g .
\end{align*}
Similarly, we have
\begin{align*}
\det(\id - g t) \det f &= \det(\id - g t) \det (\lambda g^{-1}) \\
&= (-1)^n \lambda^n t^n \det(\id - g^{-1} / t) = (-1)^n \lambda^n t^n \det(\id - f / \lambda t).
\end{align*}

\end{proof}

\begin{proof}[Proof of the main theorem]
Let $f^*: H^* (\tilde{M};\Q) \to H^* (\tilde{M};\Q)$ and $f_{\rm{rel}}^*: H^* (\tilde{M},\partial \tilde{M};\Q) \to H^* (\tilde{M},\partial \tilde{M};\Q)$ be the induced homomorphisms of $f$.
By the naturality of the cup product and Proposition~\ref{duality},
\begin{align*}
f^{*i-1}(x) \cup f_{\rm{rel}}^{*n-i}(y) &= f_{\rm{rel}}^{*n-1}(x \cup y) \\
					&= \lambda(x \cup y)
\end{align*}
is a non-degenerate pairing and satisfies the condition of Lemma~\ref{formula}. Recall that the degree $\lambda$ is with respect to the compact support cohomology $H^n_{\rm{cpt}} (\tilde{M},\partial \tilde{M})$ and the connecting homomorphism $\delta^*:H^{n-1}(\tilde{M},\partial \tilde{M}) \to H^n_{\rm{cpt}} (\tilde{M},\partial \tilde{M})$ is an isomorphism.


By using Proposition~\ref{rationality} (or Remark~\ref{rationality2}) and Lemma~\ref{formula}-(2),
\begin{align*}
\zeta_f(\frac{1}{\lambda z}) &= \prod_{i=0}^{n-1} \left[ \det(\id - f^{*i}/\lambda z) \right]^{(-1)^{i+1}} \\
	&= \prod_{i=0}^{n-1} \left[ (- \lambda z)^{-b_i} \det(\id - f_{\rm{rel}}^{*n-i-1} z) \det f^{*i} \right]^{(-1)^{i+1}} \\
	& \hspace*{8cm} (b_i: i\text{-th Betti number})\\
	&= (-\lambda z)^{\chi} \prod_{i=0}^{n-1} \left[ \det(\id - f_{\rm{rel}}^{*n-i-1} z) \right]^{(-1)^{i+1}} \times \prod_{i=0}^{n-1} \left[ \det f^{*i} \right]^{(-1)^{i+1}} \\
	&= (-\lambda z)^{\chi} \left[ \prod_{k=0}^{n-1} \left[ \det(\id - f_{\rm{rel}}^{*k} z) \right]^{(-1)^{k+1}} \right]^{(-1)^{n-1}} \times \prod_{i=0}^{n-1} \left[ \det f^{*i} \right]^{(-1)^{i+1}} \\
	&= (-\lambda z)^{\chi} \left[ \zeta^{\rm{rel}}_f (z) \right]^{(-1)^{n-1}} \times \prod_{i=0}^{n-1} \left[ \det f^{*i} \right]^{(-1)^{i+1}} .
\end{align*}
In the same way, we obtain
$$
\zeta^{\rm{rel}}_f(\frac{1}{\lambda z}) = (-\lambda z)^{\chi} \left[ \zeta_f (z) \right]^{(-1)^{n-1}} \times \prod_{i=0}^{n-1} \left[ \det f_{\rm{rel}}^{*i} \right]^{(-1)^{i+1}} .
$$

By the assumption, $n$ is odd and using Lemma~\ref{formula}-(1), we get
\begin{align*}
&\zeta^{\rm{rel}}_f(\frac{1}{\lambda z}) \zeta_f(\frac{1}{\lambda z}) \\
&= (-\lambda z)^{2\chi} \left[ \zeta^{\rm{rel}}_f (z)\zeta_f (z) \right] \times \prod_{i=0}^{n-1} \left[ \det f_{\rm{rel}}^{*i} \right]^{(-1)^{i+1}} \prod_{i=0}^{n-1} \left[ \det f^{*i} \right]^{(-1)^{i+1}} \\
&= (\lambda z)^{2\chi} \left[ \zeta^{\rm{rel}}_f (z)\zeta_f (z) \right] \times \prod_{i=0}^{n-1} \left[ \det f_{\rm{rel}}^{*i} \right]^{(-1)^{i+1}} \prod_{k=0}^{n-1} \left[ \det f^{*n-k-1} \right]^{(-1)^{k+1}} \\
&= (\lambda z)^{2\chi} \left[ \zeta^{\rm{rel}}_f (z)\zeta_f (z) \right] \times \prod_{i=0}^{n-1} \left[ \det f_{\rm{rel}}^{*i} \det f^{*n-i-1} \right]^{(-1)^{i+1}} \\
&= (\lambda z)^{2\chi} \left[ \zeta^{\rm{rel}}_f (z)\zeta_f (z) \right] \times \prod_{i=0}^{n-1} \left[ \lambda^{b_i} \right]^{(-1)^{i+1}} \\
&= (\lambda z)^{2\chi} \left[ \zeta^{\rm{rel}}_f (z)\zeta_f (z) \right] \times \lambda^{-\chi} \\
&= \lambda^\chi z^{2\chi} \left[ \zeta^{\rm{rel}}_f (z)\zeta_f (z) \right] .
\end{align*}

From Lemma~\ref{form:rel}, we complete the proof.

\end{proof}

\section{Application to knot theory}\label{sec:knot}

The Alexander polynomial is one of the most important invariants of knots. In this section, we apply our Main Theorem to the infinite cyclic covering of a knot complement and reduce the functional equation to the reciprocity of the Alexander polynomial.

Let $K \subset S^3$ be an oriented (tame) knot, $N$ a tubular neighborhood of $K$, and put $X = S^3 \setminus \Int N$. The infinite cyclic covering of $X$, i.e. the covering associated with the kernel of the abelianization homomorphism $\pi_1(X) \to H_1(X) \cong \Z$, will be denoted by $X_\infty$. By the theory of covering, the infinite cyclic group $\Pi=\langle t \rangle$ acts on $H_1(X_\infty)$, where generator $t$ is chosen to be the positive linking number with oriented knot $K$.

Identifying the integral group ring of $\Pi$ with the ring of Laurent polynomials $\Z[t^{\pm1}]=\Z[t,t^{-1}]$, $H_1(X_\infty)$ is a finitely generated $\Z[t^{\pm1}]$-module. Similarly, $H_1(X_\infty;\Q)$ is a finitely generated $\Q[t^{\pm1}]$-module.

Let $M$ be a module over a commutative ring $R$. A finite {\it presentation} for $M$ is an exact sequence
$$
F \exact{\a} E \exact{\phi} M \exact{} 0,
$$
where $E$ and $F$ are free $R$-modules with finite bases $\{e_1,\dots,e_m\}$ and $\{f_1,\dots,f_n\}$.
If $\a$ is represented by an $m \times n$ matrix, then the matrix $A$ is a {\it presentation matrix} for $M$.

By adjoining rows of zeros if necessary, we may suppose that $A$ is $m \times n$ with $m \ge n$. Then the {\it $r$-th elementary ideal} of $M$ is the ideal in $R$ generated by all the $(n-i+1) \times (n-i+1)$ minors of $A$.

\begin{Def}
The {\it $r$-th Alexander ideal} of an oriented knot $K$ is $r$-th elementary ideal of $\Z[t^{\pm1}]$-module $H_1(X_\infty)$. The $r$-th Alexander polynomial of $K$ is a generator of the smallest principal ideal of $\Z[t^{\pm1}]$ which contains the $r$-th Alexander ideal. The first Alexander polynomial is called the {\it Alexander polynomial} and is denoted by $\Delta_K(t)$.
\end{Def}

\begin{prop}[\cite{MR98f:57015}]\label{Alex=Char}
Let $K$ be a knot in $S^3$ and $t:X_\infty \to X_\infty$ be the covering transformation of $X_\infty$. Then \\
(1) $H_1(X_\infty;\Q)$ is a finite-dimensional vector space over the field $\Q$, and\\
(2) the characteristic polynomial of the linear map $t_*:H_1(X_\infty;\Q) \to H_1(X_\infty;\Q)$ is, up to a multiplication of a unit of $\Q[t^{\pm1}]$, equal to the Alexander polynomial of $K$.
\end{prop}

\begin{lem}\label{Alex=Lef}
Let $K$ be a knot in $S^3$ and $t:X_\infty \to X_\infty$ be the covering transformation of $X_\infty$. Then the Alexander polynomial $\Delta_K(z)$ and the Lefschetz zeta function satisfy the relation:
$$
\zeta_t (z) = 
\frac{1}{a_n} \frac{z^{b_1}\Delta_K(z^{-1})}{1-z},
$$
where the Alexander polynomial is normalized as $\Delta_K(z)=a_0+a_1z+\dots+a_nz^n$ with $a_0a_n\ne0$, and $b_1(=n)$ is the first Betti number of $X_\infty$.
\end{lem}

\begin{proof}
$H_i(X_\infty)=0$ for $i\ge2$ (see~\cite{MR87b:57004}).
From the Propositions \ref{rationality} and \ref{Alex=Char}, we complete the proof.
\end{proof}

From this lemma, Main Theorem implies the following well known theorem~\cite{seifert} and also \cite{MR19:53a,MR83h:58074,MR25:4526,MR15:979b}
.
\begin{cor}
The Alexander polynomial $\Delta_K(z)$ of a knot $K$ is reciprocal, i.e. $\Delta_K(z) = z^{b_1}\Delta_K(z^{-1})$.
\end{cor}

\begin{proof}
Proposition~\ref{Alex=Char}-(1) means that $X_\infty$ satisfies the requirement for the main theorem. In our condition, $\partial X_\infty \cong S^1 \times \R $. Therefore, we can see that $\zeta_{t|\partial X_\infty}(z)=1$. 
Applying Main Theorem to Lemma~\ref{Alex=Lef}, the proof is completed.
\end{proof}


\begin{rem}
From this corollary, we can rewrite the relation of Lemma~\ref{Alex=Lef} as the following:
$$
\zeta_t (z) = 
\frac{1}{\Delta_K(0)} \frac{\Delta_K(z)}{1-z}.
$$
This implies that the Alexander polynomial of a knot contains all the information about the Lefschetz numbers $\Lambda(t^n)$ of the covering transformations.
\end{rem}

\begin{ack}
I would like to express my gratitude to Professor Kazuo Masuda, Professor Hiroyuki Ochiai, Professor Hitoshi Murakami, Professor Sadayoshi Kojima, Professor Masanori Morishita and Professor Akio Kawauchi for their valuable advice.
\end{ack}

\providecommand{\bysame}{\leavevmode\hbox to3em{\hrulefill}\thinspace}
\providecommand{\MR}{\relax\ifhmode\unskip\space\fi MR }
\providecommand{\MRhref}[2]{%
  \href{http://www.ams.org/mathscinet-getitem?mr=#1}{#2}
}
\providecommand{\href}[2]{#2}

\end{document}